\documentclass[12pt,reqno]{amsart}

\usepackage{AKstyle}

\begin{document}

\title[First Non-vanishing Quadratic Twist]
{
The first non-vanishing quadratic twist of an automorphic $L$-series 
}

\author{Jeff Hoffstein
}
\thanks
{%
Hoffstein is partially supported by NSF grant DMS 0652312.
} 
\email{jhoff@math.brown.edu}
\address{Department of Mathematics,
Brown University, Providence, RI, $02912$}

\author{Alex Kontorovich} 
\thanks
{
Kontorovich is partially supported by  NSF grants DMS-0802998 and DMS-0635607, and the Ellentuck Fund at IAS} 
\email{alexk@math.brown.edu}
\address{Department of Mathematics,
Brown University, and Institute for Advanced Study, Princeton, NJ}

%, Siegel zero}%, quadratic twist} 

 \date{\today}

\begin{abstract}
Let $\pi$ be an  
automorphic representation on $\GL(r,\bA_{\Q})$ for 
$r=1, 2$, or $3$. 
Let $d$ be a fundamental discriminant and $\chi_d$ the corresponding quadratic Dirichlet character. We consider the question of the least $d$, relative to the data (level, 
weight
or 
eigenvalue) of $\pi$, such that the central value  
of the twisted $L$-series is nonzero, i.e. $L(1/2, \pi\otimes\chi_d)\neq0$. 

For example, let $N$ be the level of $\pi$. Using multiple Dirichlet series, we 
prove the nonvanishing of a central twisted $L$-value  
with $|d|\ll_\vep N^{1/2+\vep}$ for $\GL(1)$, $|d|\ll_\vep N^{1+\vep}$ for $\GL(2)$, and\break
 $|d|\ll_\vep N^{2+\vep}$ for 
$\GL(3)$, the last case assuming that a certain character is quadratic (see Theorem \ref{thm0}). We work over $\Q$ for simplicity but the method  generalizes to arbitrary number fields.

 We conjecture that in all cases there should be such a twist with $|d| \ll_\vep N^\vep$. This would follow from a Lindel\"of-type bound for a multiple Dirichlet series which does not have an Euler product, but is constructed from a Rankin-Selberg integral applied to automorphic forms which are eigenfunctions of Hecke operators.
\end{abstract}

\maketitle

\tableofcontents

%\newpage

\section{Introduction}

Much work has been devoted to the problem of bounding the least prime in an arithmetic progression. The Grand Riemann Hypothesis %implies 
predicts
that for $(a,q)=1$ and any $\vep>0$,
\be\label{eq:PNTAP}
\sum_{p\equiv a(q)\atop p<N}\log p = {N\over \varphi(q)} + O_{\vep}\bigg( N^{1/2+\vep}\bigg),
\qquad
\text{ as }\qquad N\to\infty.
\ee
%with the implied constant depending only on $\vep$. 
As $\varphi(q)=q^{1+o(1)}$, %\eqref{eq:PNTAP} 
the above
implies that the main term dominates the error as soon as $N\gg_{\vep}q^{2+\vep}$. Hence the left hand side of \eqref{eq:PNTAP} is non-zero, confirming the existence of a prime 
\be\label{eq:pllq2}
p\ll_{\vep} q^{2+\vep}
\ee 
with $p\equiv a(q)$.

Unconditionally, the error term in \eqref{eq:PNTAP} is not much better than a power of log savings, the consequence being an exponential rather than polynomial bound in \eqref{eq:pllq2}. Linnik \cite{Linnik1944} was the first to show that the problem of least prime in an arithemetic progression is not subordinate to progress towards the ``first moment'' in \eqref{eq:PNTAP}, proving a polynomial bound (the quality of which has since been vastly improved).
\\

In this paper we consider the following similar problem. Let $\pi$ be an automorphic 
%form 
representation on $\GL(r,\bA_{\Q})$, $r=1, 2$ or $3$, and let $\cD$ denote the set of fundamental discriminants.
For $d\in\cD$, let $\chi_d$ be the corresponding quadratic Dirichlet character of modulus $d$.
A great deal of attention has been paid in recent years to the question of the existence and  abundance of % fundamental discriminants 
$d\in\cD$ such that 
the central value of
the standard $L$-function attached to $\pi$ twisted by $\chi_{d}$,
$L(1/2,\pi\otimes\chi_d),$ does not vanish.
%, that is, such that the  $L$-series of $f$, twisted by $\chi_d$, does not vanish at the center of the critical strip.  
%In this paper, w
We pose the following more refined problem. %question:

\begin{question}%[Q]
\label{Q}
{\rm
If such a $d$ exists, then what is the least value of  $|d|$, relative to the data of $\pi$ (such as its level $N$ %, weight $k$, 
or eigenvalue $\gl$), for which the twisted $L$-series does not vanish at the center?%
}
\end{question}

\begin{rmk}
%We point out first that
%We immediately dismiss
Though one might {\it a priori} assume that
 the ``analytic conductor'' of Iwaniec-Sarnak %is {\bf not} 
%as
is
 a suitable %appropriate 
 measure of the ``complexity'' of $\pi$ in this 
problem,
 %situation,
 the following examples suggest %otherwise
  %-- 
  in fact that
  the dependence on the archimedean place plays a substantially %vastly 
 different 
 role from the %level
finite ramification. %! %finite places
% Here is some %immediate 
% ``trivial''
% progress towards Question \ref{Q}, illustrating the distinction.
 
 \begin{enumerate}
\item
%It should be noted that for  
For
%$r=2$ and 
$\pi$ a discrete series representation of $\GL(2)$ corresponding to a holomorphic modular form of level $N$, by simply combining  Waldspurger's theorem \cite{Waldspurger1981} with Riemann-Roch, one can  show the existence of a $|d|\ll_{\vep} N^{1+\vep}$ %(no epsilon!)
 such that the central $L$-value is non-zero. 
 
 \item 
% For  $r=2$, %as before, 
Similarly,
for $\pi$ a tempered representation of $\GL(2)$ having Casimir eigenvalue $\gl,$%
\footnote{%
Classically, this corresponds to either a Maass form of eigenvalue 
$\gl=\frac14+t^{2}$, $t>0$, or a holomorphic form of even weight $k$ with
$\gl=\frac k2(1-\frac k2)$.
}
%a  principal (resp. discrete) series representation of $\GL(2)$ %spherical having a $K$-%fixed vector 
%corresponding to a Maass (resp. holomorphic) form $f$ 
%of 
%eigenvalue 
one can
%again
 %show via
 apply
  Waldspurger's theorem and arguments dating back to Maass%
\footnote{P. Sarnak pointed us to % the thesis of 
J. Huntley's thesis, where these ideas are vastly generalized.}
%These imply 
to show 
the existence of
%that an automorphic form of even weight $k$  has 
a nonvanishing central twist $\chi_d$ with 
%$|d|\ll k$, and that a Maass form of eigenvalue $\gl=\frac14+t^2$ has a nonvanishing twist with 
$|d|\ll %t\approx
\gl^{1/2}$. %Similarly, 
 \end{enumerate}
%It is interesting to note that while 
As
the analytic conductor in the above examples is roughly $%C(\pi)
%\approx 
N\gl$,
%
%in the level aspect is $N$
%,   in the weight aspect %it 
%is $k^2$ 
%and  is $\gl$ for the eigenvalue
% . Despite the fact that ``convexity'' arguments apply identically to both aspects, 
 %the Multiple Dirichlet Series technology sees the difference, and knows that 
 %it is the dimension of the space (being $\sim N$ or $\sim k$), not the conductor, which
it is clear that  one 
%must 
should
separate  the level and eigenvalue aspects
in %studying Question \ref{Q}
this problem. That said, see the caveat in Remark \ref{rmk:Conductor}. 
%  plays the role of being the convex bound.
\end{rmk}

%\begin{rmk}
The  examples above can be considered ``convexity'' bounds towards Question \ref{Q}, for reasons which 
shall become clear, see \S\ref{sec:conv}. %Remark \ref{rmk:conv}.
%\end{rmk}
%
 In this paper, we demonstrate the convexity bound for the central $L$-value of 
$\pi$ on $\GL(r)$ with $r=1,2,$ and $3$.

\subsection{Statements of the Main Results}
\

%Our assumptions on $\pi$ are as follows. 
We shall really only work with the standard $L$-function attached to quadratic twists of $\pi$,
% and its quadratic twists, 
and not $\pi$ itself.
 %(recall that twists by all Dirichlet characters suffice to determine $\pi$ by the Converse Theorem). 
 To this end, we make the following
%The standard $L$-function associated to $\pi$ is
\begin{Def}
By an automorphic $L$-series $L(s,\pi)$ on $\GL(r)$, we mean the following.
Assume that the series
$$
L(s,\pi) = \sum_{n\ge1}^\infty \frac{c(n)}{n^s}
$$
converges absolutely for $\Re(s)$ sufficiently large,
has Euler product
\be\label{eq:EulerProd}
L(s,\pi) = \prod_p  \prod_{j=1}^r\big( 1 - \alpha_p^{(j)} p^{-s}\big)^{-1} ,
\ee
and analytic continuation with functional equation 
\be\label{functeqplain}
\gL(s,\pi):=N^{s/2}G_\pi(s) L(s,\pi)  = 
 \epsilon_\pi \gL(1-s,\tilde \pi) %N^{(1-s)/2}G_f(1-s) L(1-s,\tilde f)
 .
\ee
Here $\tilde \pi$ is the 
contragradient of $\pi$, $|\epsilon_\pi|=1$ is the root number, the integer $N\ge1$ is the level, % analytic conductor, 
and
\be\label{eq:GamFac}
G_\pi(s) = \pi^{-rs/2}\prod_{j=1}^r\G\left( {s+\gk_j\over 2} \right)
\ee
is a product of archimedean gamma factors with $\Re(\gk_{j})\ge0$. 
We define the ``archimedean conductor'' of $\pi$ to be:
\be\label{eq:glDef}
\fq%_{\infty} 
:= \prod_{j=1}^r\left( 3+|\gk_j| \right).
\ee
For positive square-free $d$, the twisted $L$-series has Euler product
\be\label{Ldef}
L(s,\pi\otimes\chi_d)  = \sum_n {c(n)\chi_d(n)\over n^s} =  \prod_p  \prod_{j=1}^r\big( 1 - \chi_d(p) \alpha_p^{(j)} p^{-s}\big)^{-1}
,
\ee
and %if $(N,d)=1$ then 
its functional equation is given by
\bea
\label{functeq}
\gL(s,\pi\otimes \chi_d)&:=&(
%N|D|^r
\fc(\pi\otimes \chi_{d})
)^{s/2}G_{d,\pi}(s) L(s,\pi\otimes\chi_d)  
\\
\nonumber
&=&
 \epsilon_\pi \psi(d) \gL(1-s,\tilde \pi\otimes\chi_d) 
 %(N|D|^r)^{(1-s)/2}G_{d,f}(1-s) L(1-s,\tilde f\otimes\chi_d)
 .
\eea
Here
$G_{d,\pi}$ is a product of gamma factors depending only on $\pi$ and the sign of $d$, %and 
$\psi$ is a character modulo $N$,
and $\fc(\pi\otimes\chi_{d})$ is the conductor of the twisted automorphic representation $\pi\otimes\chi_{d}$. %, not necessarily primitive.
For example, if $(N,d)=1$ then
\vskip-.23in
\be\label{eq:unramCond}
\fc(\pi\otimes\chi_{d})=N|D|^{r},
\ee 
where
 $D = 4d$ or $D=d$ is the conductor of $\chi_d$.
%
%
%If 
For the reader's convenience, we record %classify 
the possible values of $\fc(\pi\otimes\chi_{d})$ in the case
$(N,d)\ne1$
%, the factor $N|D|^r$ is replaced by the conductor of 
%$\pi\otimes \chi_d$.   See 
in
\S\ref{sec:Din}. % for a detailed discussion.
\end{Def}

We first state our main result in  the level aspect:

\begin{theorem}\label{thm0}
Let $L(s,\pi)$ be an automorphic  $L$-series on $\GL(r)$ 
of level $N$ and degree
$r=1$, $2$ or $3$.  Suppose that the root number of   $L(s,\pi)$ is \emph{not} equal to $-1$ (and hence there exists a non-vanishing quadratic twist, see \S\ref{sec:Siegel}).  Suppose also in the case %that if 
$r=3$ that %then 
$\psi$ is trivial or quadratic.\footnote{%
This is but a simplifying assumption;
for the most general statement, see \S\ref{sec:proof}.
}
Then 
\begin{itemize}
\item
for  $r=1$, there exists some $|d|\ll_\vep N^{1/2+\vep}$,
\item
for  $r=2$, there exists some $|d|\ll_\vep N^{1+\vep}$, and
\item
for  $r=3$% and $\pi$ self-dual
, there exists some $|d|\ll_\vep N^{2+\vep}$,
\end{itemize}
such that  $L(1/2,\pi\otimes\chi_d) \ne 0$. 
%If one restricts to an arithmetic progression of $d$ for which the root number is \emph{always} $-1$, (so all central values are forced to vanish), the same result holds for  $L'(1/2,\pi\otimes\chi_d) \ne 0$.
\end{theorem}

In the eigenvalue aspect (which only makes sense over $\Q$ for degree $r=2$ and $r=3$), we have:

\begin{thm}\label{thm1}
Let $L(s,\pi)$ be an automorphic $L$-series on $\GL(r)$ of  archimedean conductor $
\fq%_{\infty}
$ as in \eqref{eq:glDef} and degree $r=2$ or $3$.
 Suppose that the root number of   $L(s,\pi)$ is not equal to $-1$.
Then 
\begin{itemize}
\item
for  $r=2$, there exists some $|d|\ll_\vep \fq%_\infty)
^{1/2+\vep}$, and
\item
for  $r=3$% and $\pi$ self-dual
, there exists some $|d|\ll_\vep \fq%_\infty)
^{1+\vep}$,
\end{itemize}
such that  $L(1/2,\pi\otimes\chi_d) \ne 0$.
\end{thm}

\begin{rmk}
The archimedean conductor $\fq%_\infty)
$ should not be confused with the Casimir eigenvalues of $\pi$.
Consider a principal series representation $\pi$   on $\GL(3)$ corresponding to a 
 Maass form of type $\nu=(1/3+it_1,1/3+it_2)$, with $t_j\asymp T$. 
 The eigenvalue of the Laplacian %Casimir operator %$\cC$ 
  is then
 $$
\gl%_{\cC}
=
1 + 3 t_1^2 + 3 t_1 t_2  + 3 t_2^2
 \asymp T^2
 .
 $$ 
The Gamma factors of $\pi$ are (cf. \cite[Theorem 6.5.15]{Goldfeld2006})
$$
G_\pi(s):=
\G\left(
s+2it_1+it_2
\over
2
\right)
\G\left(
s+it_1-it_2
\over
2
\right)
\G\left(
s-it_1-2it_2
\over
2
\right)
,
$$
so generically the archimedean conductor is
$$
\fq%_\infty)
\asymp
T^3
\asymp
\gl%_{\cC}
^{2/3}
.
$$
Then Theorem \ref{thm1}  exhibits a nonvanishing twist $\chi_{d}$ with 
$$
|d|\ll_\vep \gl%_{\cC}
^{3/2+\vep}.
$$   
On the other hand, on $\GL(2)$, $\fq%_\infty)
\asymp\gl%_{\cC}
$.
\end{rmk}

\subsection{Outline of the Proof}

\subsubsection{The ``Moment'' Method}\label{sec:Moment}\

Let $L(s,\pi)$ be an automorphic $L$-series on $\GL(r)$, $r=1, 2,$ or $3$.
As in the case of primes in progressions, one can try to 
compute the first moment:
\be\label{eq:MomentL}
\sum_{d\in\cD\atop|d|<X}L(1/2,\pi\otimes\chi_{d})=M_{\pi}(X) + E_{\pi}(X).
\ee
%and \pi
%(Usual
For $X$ large enough that the main term dominates the error, the above formula will produce a non-vanishing central twist.
%$M_{\pi}(X)>|E_{\pi}(X)|$
In practice, 
$$
M_{\pi}(X)=X^{1+o(1)}
,
$$ 
and in order to prove, say Theorem \ref{thm0} in the level $N$ aspect, one needs to bound the error term $E_{\pi}(X)$ by terms of the form 
$$
X^{1-\ga} N^{\ga\cdot \gt%_{r}
}
,
$$ 
for some $0<\ga<1$ with  $\gt%_{1}
=1/2$, $\gt%_{2}
=1$, or $\gt%_{3}
=2$ corresponding to $\GL(1)$, $\GL(2)$, or $\GL(3)$, respectively.
Even with smooth weights, 
unconditional
moments with this quality of error seem difficult to achieve with existing methods, especially on higher rank groups such as $\GL(3)$ (see Remark \ref{rmk:Moment}). As in %the case of 
Linnik's problem, %one can
we will
 establish first non-vanishing results {\it without}
making progress towards \eqref{eq:MomentL}. %, and this will be our approach.

\subsubsection{The ``Multiple Dirichlet Series''  Method}\

Instead of computing the moment, we employ the theory of double Dirichlet series. 
Consider the following Dirichlet series, whose coefficients are themselves twisted $L$-functions with some carefully chosen weights:
\be\label{eq:Zdef}
Z(s,w):=\sum_{d}{L(s,\pi\otimes\chi_{d})\ P(s,\pi,d)\over d^{w}}.
\ee
The series thus defined converges for $\Re(s),\Re(w)$ sufficiently large. As has been detailed in many places (e.g. \cite{DiaconuGoldfeldHoffstein2003, BumpFriedbergHoffstein2004} etc.), $Z(s,w)$ has meromorphic continuation to all $(s,w)\in\C^{2}$ with explicitly understood polar divisors, and satisfies a finite group of functional equations, including the transformation 
\be\label{eq:swTo1s1w}
(s,w)\mapsto (1-s,1-w).
\ee
Specializing to $s=1/2$, one obtains a functional equation of the form
\be\label{eq:MDSFE}
G(w)N^{\gt w}Z(1/2,w) \approx \tilde G(1-w) N^{\gt (1-w)}\tilde Z(1/2,1-w),
\ee
where $G$ and $\tilde G$ are archimedean (Gamma) factors and $\tilde Z$ is constructed in a similar way as $Z$. 
(The true functional equation is actually a linear combination of terms like $\tilde Z$, \`a la the ``scattering matrix'' in the functional equation of an Eisenstein series, see 
e.g. equations \eqref{funct1} -- \eqref{funct3}.) 
Moreover, $Z(1/2,w)$ has a pole at $w=1$ (and possibly at $w=3/4$ %and $w=1/4$ 
on $\GL(3)$). By a familiar Tauberian argument resembling an approximate functional equation, we can thus write the residue at $w=1$ as a finite sum of coefficients of $Z(1/2, w)$ (which are of course the sought-after central $L$-values), where the length of the sum is the square root of the ``conductor'' in the functional equation \eqref{eq:MDSFE} for the double Dirichlet series:
\be\label{eq:MDSAFE}
\text{Residue}
\approx
 \sum_{d}{L(1/2,\pi\otimes\chi_{d})\ P(1/2,\pi,d)\over d^{1/2}}\ V\left({|d|\over N^{\gt}}\right)
 +
 \text{similar}
% +
% \sum_{d}{L(1/2,\tilde\pi\otimes\chi_{d}) \tilde P(1/2,\tilde \pi,d)\over d^{1/2}}\tilde V\left({d\over N^{\gt/2}}\right)
,
\ee
where $V$ is supported in $[1,2]$, say.
Then one immediately arrives at a contradiction if all $L$-values vanish with $|d|\ll N^{\gt+\vep}$. 
The same argument applies to the archimedean aspect.% 
\footnote{%
Also, functional equations of the type \eqref{eq:MDSFE} have  been worked out in complete detail over number fields \cite{BumpFriedbergHoffstein2004}, so our approach %immediately 
extends to this setting. For ease of exposition, we will restrict ourselves to $\Q$. 
}

\begin{rmk}\label{rmk:Conductor}
In fact one {\it can} combine the level and eigenvalue Theorems \ref{thm0} and \ref{thm1} into a uniform statement, but this  involves knowing the ``conductor'' for the double Dirichlet series attached to $\pi$, and not just $\pi$ itself. It is in general difficult to predict {\it a priori} the exact shape of the functional equation in \eqref{eq:MDSFE} without following through a sequence of functional equations as in \S\ref{sec:FEheur} to reach the transformation \eqref{eq:swTo1s1w}. So the appropriate ``conductor'' cannot elementarily be read off from the functional equation \eqref{functeqplain}. See also Remark \ref{rmk:condGuess}.
\end{rmk}

\subsection{Subconvexity}\label{sec:conv}
\

Note that the approximate functional equation \eqref{eq:MDSAFE} for the double Dirichlet series is morally equivalent to a ``convexity'' bound for the series $Z(1/2,w)$ at the central point $w=1/2$.
Indeed, the Tauberian argument alluded to before is to consider essentially
\bea
\label{eq:contour}
&&
\hskip-1in
{1\over 2\pi i}\int_{(2)} G(w+1/2) %N^{\gt (w+1/2)} 
Z(1/2,w+1/2) X^{w}dw
\\
\nonumber
&&
=
\sum_{d} {L(1/2,\pi\otimes\chi_{d})\ P(1/2,\pi,d)\over d^{1/2}} \ V\left({|d|\over X
%\over N^{\gt}
}\right)
.
\eea
Pull the contour past the pole at $w=1/2$ all the way to the line $\Re(w)=-1$, say, after which applying the functional equation and taking $X=N^{\gt}$ recovers \eqref{eq:MDSAFE}.

If one were to attempt an improvement via these methods on the exponents in  Theorems \ref{thm0} and \ref{thm1},
a key ingredient would be a ``subconvex'' %an
 bound for  $Z(1/2,w)$ at %the center of its critical strip, i.e. at 
 $w=1/2$.   
% 
 %The upper bound that corresponds to the first non-vanishing ranges described above is  obtained by a Phragm\'en-Lindel\"of convexity argument.   
 Any improvement on the convexity bound would lead to a corresponding improvement of these results, and a full Lindel\"of-type bound would lead to the existence of a non-vanishing twist  with 
$$
|d|
\ll_\vep (\fq N)^\vep.
$$ 

 The interesting point is that the $L$-series $Z(1/2,w)$ does {\bf not} have an Euler product.   
  Of course one does not expect a Lindel\"of-type bound to be true in general for $L$-series without an Euler product.  See e.g., \cite{ConreyGhosh2006} where a counterexample is constructed.  However, it does not seem unreasonable to conjecture a Lindel\"of type bound for an $L$-series without an Euler product when that $L$-series  is constructed from a Rankin-Selberg integral %(granted, with respect to an Eisenstein series on a metaplectic cover) 
  applied to one or more automorphic forms that are themselves eigenfunctions of the relevant Hecke operators.  Indeed, we conjecture  that the  double Dirichlet series $Z(s,w)$ satisfies 
\be\label{eq:conj}
Z(1/2,1/2+it)\ll_\vep (\fq N)^{\vep}(1+|t|)^{A},
\ee
for some $A>0$.
Were this to be the case, one could just pull the contour in \eqref{eq:contour} to the line $\Re(w)=0$ (still collecting the residue at $w=1/2$) and estimate away the remaining integral. Since the variable $X$ is free, one can choose it to be as %large
small
 as $N^{%\gt-
 \vep}$, making the sum on the right hand size of \eqref{eq:contour} have negligible length,% 
\footnote{%
Note added in print: 
In  recent work, Blomer \cite{Blomer2009} has succeeded on $\GL(1)$ and level $N=1$ in proving a subconvex estimate for a double Dirichlet series in the $t$-aspect. Unfortunately this  is the only aspect which does not give applications towards our questions. Of course, it does give evidence that more progress is within reach.
}
cf. Remark \ref{rmk:SubConv}.

\begin{rmk}\label{rmk:Moment}
Unconditionally, one has some polynomial bound in \eqref{eq:conj}, see   \eqref{eq:convR1} -- \eqref{eq:convR3}. One can start with the left hand side of \eqref{eq:contour}, except without the Gamma factors, pull the line to $\Re(w)=0$, and estimate the error there after extracting the residue. Note that this {\it does not} recover the same result as Theorems \ref{thm0} and \ref{thm1}! In fact, this approach is much closer to that of the ``moment'' approach described in \S\ref{sec:Moment}. %, and leads to correspondingly weaker results.
\end{rmk}

\subsection{Degree $r\ge4$}\label{r4}
\

%{\rm
On $\GL(r)$ with $r\ge4$, the group of functional equations is no longer a finite Weyl group, but is an infinite Coxeter group, see Remark \ref{rmk:gl4}. %The c
Current technology is  incapable in this case of obtaining the analytic continuation of $Z(s,w)$ beyond the critical point $(s,w)=(1/2,1)$, the sole exception being the recent work by Bucur and Diaconu \cite{BucurDiaconu2008}  in the  function field analogue. Moments for quadratic twists of generic $\pi$ on $\GL(4)$ and higher are also presently unavailable. 
 
 In particular, one cannot yet 
answer the following enticing question.
Given two automorphic forms $\pi$ and  $\pi'$ on $\GL(2)$, each with a positive sign in their functional equation, does there exist a quadratic twist $\chi_d$ such that the two twisted $L$-series simultaneously do not vanish at the center of the critical strip, i.e. %such that 
$L(1/2,\pi\otimes\chi_d) L(1/2,\pi'\otimes\chi_d)\ne 0$?  
Similarly, one cannot yet  
 obtain  the second moment of an automorphic form $\pi$ on $\GL(2)$ twisted by quadratic characters, i.e. an asymptotic formula for 
 $$
 \sum_{d<X} L(1/2,\pi\otimes\chi_d)^2,\quad\quad\text{ as $X\to\infty$.}
 \footnote{
Note added in print: in the recent work \cite{SoundYoung2009},  Soundararajan and Young give unconditional {\it lower} bounds  which match the conjectured main term!
 }
 $$
 
%}
%\end{remark}

%\begin{remark}
%{\rm

\subsection{Moments of Half-integral Weight Forms}
\

By the Shimura correspondence, the questions raised above for $\GL(2)$ are %not un
related to questions about twisted moments of half-integral weight forms. Again, the $L$-series attached to a half-integral weight form $\tilde f$ does not have an Euler product, yet it seems likely that if the integral weight Shimura correspondent $f$ is an eigenfunction of the Hecke operators, then $L(s,\tilde f)$ should satisfy a Lindel\"of type bound at the center of its critical strip.
In joint work with Gautam Chinta  \cite{ChintaHoffsteinKontorovich2010}, we have observed that, contrary to the integral weight situation, if one forms the multiple Dirichlet series
$$
\tilde Z(s_1,s_2,w)\approx \sum_d {L(s_1,\tilde f\otimes\chi_d)L(s_2,\tilde f\otimes\chi_d) \over d^w},
$$
then its group of functional equations is isomorphic to the Weyl group associated to the Dynkin diagram $A_5$, which is finite!  Thus we are able to
obtain first and second moments for half-integral weight forms twisted by quadratic characters, i.e. asymptotics as $X\to\infty$ for
$$
\sum_{d<X}L(1/2,\tilde f\otimes\chi_d)
\quad\text{ and }\quad
\sum_{d<X}L(1/2,\tilde f\otimes\chi_d)^2. 
$$
As the first pole of $\tilde Z(1/2,1/2,w)$ appears at $w=1$, the second moment is asymptotic to 
$
X\, P(\log X)
,
$
where $P$ is some polynomial. This gives further evidence of the truth of a Lindel\"of type bound, even for certain ``arithmetic'' $L$-functions without Euler products.
%}
%\end{remark}

\subsection{Simultaneous Non-vanishing Twists}
\

Choosing the representation $\pi$ %on $\GL(2)$ or $\GL(3)$ 
in a particular way, such as 
$$
\pi=\chi_{N_{1}}\boxplus\chi_{N_{2}}\boxplus\chi_{N_{3}}\qquad\text{ or }\qquad\pi=\pi_{1}\boxplus\chi_{N_{1}}
$$ 
for characters $\chi_{N_{j}}$ and  $\pi_{1}$ on $\GL(2)$,
Theorem \ref{thm0} has the following immediate corollary on simultaneously non-vanishing twists.
\begin{corollary}
Let $L(s,\chi_{N_1}), L(s,\chi_{N_2}),L(s,\chi_{N_3})$ be three Dirichlet $L$-series
%, not necessarily quadratic, 
with conductors $N_1,N_2, N_3$.    Let $L(s,\pi)$ be an automorphic $L$-series on $\GL(2)$ of level $N$.
%Assume that $L(s, \chi_{N_i} \chi_{N_j})$ have no Siegel zeros for $1\le i<j\le 3$. 
Then
\begin{enumerate}
\item
 there exists 
 %a $d$, with  
  $|d| \ll  (N_1 N_2)^{1+ \epsilon}$ with 
$$L(1/2,\chi_d\chi_{N_1}) L(1/2,\chi_d\chi_{N_2}) \ne 0,$$
\item
there exists  $|d| \ll (N_1 N_2N_3)^{2+ \epsilon}$ such that 
$$L(1/2,\chi_d\chi_{N_1}) L(1/2,\chi_d\chi_{N_2}) L(1/2,\chi_d\chi_{N_3})\ne 0,$$ 
\item there exists %a 
$|d| \ll  (N_1 N)^{2+ \epsilon}$ such that 
$$L(1/2,\chi_d\chi_{N_1}) L(1/2,\pi\otimes\chi_d) \ne 0.$$
\end{enumerate}
%Dropping the Siegel zero assumption increases the exponents in $(1)$ and $(2)$ by $1/2$. (The exponent in $(3)$ is not affected.)
\end{corollary}

%These are merely curiosities for which we have no particular applications.

\subsection{Outline of the Paper}
\

In \S\ref{sec:prelim}, we present the heuristic derivation of  the functional equations for the double Dirichlet series $Z(s,w)$, and then state them rigorously in \S\ref{sec:FEs}. These are well-known to the experts, but our application requires slightly more refined information%, which we give in \S\ref{sec:refine}
; the level aspect is the only case which causes difficulty. Equipped with this data, we prove the main theorems in \S\ref{sec:proof}.

%\newpage

\subsection*{Acknowledgements}
The authors wish to thank Adrian Diaconu, Peter Sarnak, and Leo Goldmakher  for many comments and corrections to an earlier draft. Many thanks also to 
 Dinakar Ramakrishnan for   
 the analysis in \S\ref{sec:Din}.
 We are especially grateful to the referees for their careful reading of the manuscript, and many suggestions and improvements to the text.

\section{Preliminaries}\label{sec:prelim}

\subsection{The Heuristic Argument}\label{sec:Heur}
\

Before presenting the (quite technical) details of the functional equation leading to \eqref{eq:MDSFE}, we give a heuristic argument, % for the uninitiated, 
focusing on  the level aspect.   It will contain some very imprecise statements regarding functional equations but should nevertheless be a useful reference guide for the actual proofs.   We pretend throughout this section, for clarity of exposition, that all numbers are positive and congruent to 1 modulo 4, and that quadratic reciprocity is perfect.   We will also suppress the weights $P(s,\pi,d)$; they appear in every equation and contribute little to  the exposition.
 
Let $\pi$ be an automorphic representation on $\GL(r,\bA_{\Q})$, with $r=1, 2,$ or $3$ and Fourier coefficients $c(n)$.
Consider the following double Dirichlet series:
\be\label{step1}
Z%_1
(s,w) = \sum_d \frac{L(s,\pi\otimes\chi_d)}{d^w}.
\ee
%When $d$ is square free, $L(s,\pi\otimes\chi_d)$ is the twisted $L$-series of  \eqref{Ldef}.   When $d$ is not square free, $L(s,\pi\otimes\chi_d)$ is the modified $L$ series given by \eqref{Pdef} (taking $M=a_1 = l_1 =1$ in this simplified setting).
% 
Very roughly, inserting \eqref{Ldef} into \eqref{step1}, $Z%_1
(s,w)$ is represented by the double Dirichlet series
$$
\sum_{d,n}\frac{c(n)\chi_d(n)}{d^wn^s}.
$$
This suggests that if quadratic reciprocity held perfectly, that is $\chi_d(n) = \chi_n(d)$, then we could rewrite this as
\be\label{step2}
Z%_1
(s,w) = \sum_n \frac{L(w,\chi_n)c(n)}{n^s},
\ee
and in fact this interchange can be made rigorous, cf. \S\ref{sec:interch}.

%Dualizing
Applying the functional equation to 
the numerator of \eqref{step2} and suppressing gamma factors, 
we see that there is a functional equation sending
\be\label{step5}
Z%_1
(s,w) \rightarrow Z%_1
( s+w -1/2,1-w).
\ee
On the other hand, if we apply  \eqref{functeq} 
 to the numerator of \eqref{step1}, we find that there is a functional equation sending
\be\label{step3}
Z%_1
(s,w) \rightarrow N^{1/2 - s}\tilde Z%_2
(1-s, w + rs -r/2),
\ee
where
\be\label{step4}
\tilde Z%_2
(s,w)= \sum \frac{L(s,\tilde \pi\otimes\chi_d)\psi(d)}{d^w}
=\sum_{n,d} \frac{\chi_n(d)\psi(d)\tilde c(n)}{n^s d^w}
.
\ee
Reversing orders of summation and collecting terms in \eqref{step4}, we find that 
$$
\tilde Z%_2
(s,w)= \sum_n \frac{L(w,\chi_n\psi)\tilde c(n)}{n^s}.
$$
%Dualizing
\begin{rmk}\label{rmk:condPsi}
%Now,
In the above,
 the conductor of $\psi$ could be any divisor of $N$.   In this heuristic we will assume for simplicity that $\psi^2 =1$ and that the conductor equals $N$.
\end{rmk}
Applying the functional equation to 
this numerator, we see that there is 
a transformation %functional equation
\be\label{step6}
\tilde Z%_2
(s,w) \rightarrow N^{1/2 - w}\tilde Z%_2
( s+w -1/2,1-w).
\ee

\begin{rmk}
Note that if $\psi$ is complex, then a new character and a Gauss sum could be introduced in this functional equation.   This is why the result for $\psi$ complex is slightly worse than the result for $\psi$ real in the case $r=3$, cf. \eqref{funct3new}.  
\end{rmk}

Similarly, \eqref{step3} can be used in reverse to give
\be\label{step7}
\tilde Z%_2
(s,w)
\to
N^{1/2-s}
Z%_1
(1-s, w + rs -r/2) 
.
\ee

\subsubsection{Iterating %the Local
 Functional Equations% to get a Global One
 }\label{sec:FEheur}
 \
 
We now apply these functional equations in sequence.   If the degree $r = 1$, we apply in succession
 \eqref{step5}, \eqref{step3}, and \eqref{step6}, 
 obtaining
\bea
  \nonumber
Z%_1
(s,w)& \rightarrow &Z%_1
(s+w-1/2, 1-w) \rightarrow N^{1-s-w}\tilde Z%_2
(3/2 - s-w,s)
  \\
  \label{gl1ref}
&  \rightarrow& N^{3/2-2s-w}\tilde Z%_2
(1-w,1-s).
\eea

\begin{remark}\label{rmkGL1}
{\rm
On $\GL(1)$ there is an extra symmetry in \eqref{step4}, namely $\tilde Z%_2
(s,w)\approx \tilde Z%_2
(w,s)$, coming from the relation $\psi(d)= c(d)$.
}
\end{remark}

If $r = 2$ we apply in succession \eqref{step5}, \eqref{step3}, \eqref{step6} and \eqref{step7}, obtaining
\beann
Z%_1
(s,w) &\rightarrow &Z%_1
(s+w-1/2, 1-w) \rightarrow N^{1-s-w}\tilde Z%_2
(3/2 - s-w,w+2s -1)
  \\
&  \rightarrow & N^{5/2-3s-2w}\tilde Z%_2
(s,2-2s-w) \\
& \rightarrow &  N^{3-4s-2w}Z%_1
(1-s,1-w).
\eeann

If $r = 3$ we apply in succession \eqref{step5}, \eqref{step3}, \eqref{step6}, \eqref{step7},  \eqref{step5} and \eqref{step3}, obtaining
\beann
Z%_1
(s,w)& \rightarrow & Z%_1
(s+w-1/2, 1-w) \\
&\rightarrow& N^{1-s-w}\tilde Z%_2
(3/2 - s-w,3s+2w -2)   \\
& \rightarrow & N^{7/2-4s-3w}\tilde Z%_2
(2s+w -1,3-3s-2w) \\
&\rightarrow&  N^{5-6s-4w}Z%_1
(2-2s-w,w+3s-3/2)\\
& \rightarrow&  N^{5-6s-4w}Z%_1
(s,5/2-3s-w) \\ 
&\rightarrow&  N^{11/2-7s-4w}\tilde Z%_2
(1-s,1-w).
\eeann

\begin{remark}\label{rmk:gl4}
%{\rm
In each of the cases above, 
%enough 
finitely-many
iterations of the functional equations will return us to $Z%_1
(s,w)$. For degree $r\ge4$, one can cycle the transformations ad infinitum, never arriving at the desired argument $Z(1-s,1-w)$.
%}
\end{remark}

When the functional equations above are applied to  $Z%_1
(1/2,w)$ we find the following relations hold (in each case below, ``$\rightarrow$'' indicates that only the archimedean contributions are suppressed).
When $r=1$:
$$
N^{w/2} Z%_1
(1/2,w) \rightarrow N^{(1- w)/2}\tilde Z%_2
(1/2,1-w),
$$
when $r=2$:
$$
N^{w}Z%_1
(1/2,w) \rightarrow N^{1 - w}Z%_1
(1/2,1-w),
$$
and when $r=3$:
$$
N^{2w}Z%_1
(1/2,w) \rightarrow N^{2 (1- w)}\tilde Z%_2
(1/2,1-w)
.
$$
These are exactly the relations corresponding to \eqref{eq:MDSFE}.
\begin{rmk}\label{rmk:condGuess}
As noted in Remark \ref{rmk:Conductor}, it is a bit delicate to determine the exact form of the functional equation for $Z%_{1}
(1/2,w)$, and hence its ``analytic conductor'', given the initial data of $\pi$. But one can use the above heuristic as a template to predict the outcome.
\end{rmk}

%We now give the precise formulation of the functional equation for our double Dirichlet series. For this, we assume familiarity with the many sources 

\subsection{Ramified Conductor}\label{sec:Din}\

In this section, we give precise details for the conductor of $\pi\otimes \chi_{d}$ appearing in \eqref{functeq} in the case of ramified twist.
We are indebted to Dinakar Ramakrishnan for providing us with the following case by case analysis. % of the conductor in the case of ramified twist.
 
Let
$\pi$ %is 
be
an irreducible admissible representation of $\GL_r(\Q_p)$ and $\chi$ a character of $\Q_p^\times$. (The global conductor  is a product of such local conductors, all but finitely many of which are unity.)
Recall that  $\fc$ denotes the conductor of a representation.
The case when the conductors of $\pi$ and $\chi$ are relatively prime is trivial, and one has $\fc(\pi\otimes\chi)=\fc(\pi)\fc(\chi)^r$, as in \eqref{eq:unramCond}. % in this case.
 
Let us assume from now on that $\fc(\pi)=p^a$ and $\fc(\chi) =p^b$. The representation $\pi\otimes\chi$ is hence ramified, and one can appeal to 
Tadic's classification \cite{Tadic1986} of such. 
 There is a partition $r=r_1+r_2+...+r_m $, and discrete series representations $\pi_j, 1\leq j\leq m,$ of $\GL_{r_j}(\Q_p)$ such that $\pi$ is parabolically induced from the representation $\pi_1\times\pi_2 \times...\times \pi_m$ of the parabolic $P$ attached to the partition. One has
$$
\fc(\pi) = \prod_{j\leq m} \fc(\pi_j)
$$ 
and
$$
\fc(\pi\otimes\chi)  =  \prod_{j \leq m} \fc(\pi_j\otimes\chi).
$$ 
Thus it suffices to understand
 $\fc(\eta\otimes\chi)$ for a discrete series representation $ \eta$ of $\GL_t(\Q_p)$ and a character $\chi$.
For simplicity, assume that $t < p$.
 
Suppose $t=1$.  Then  $\eta$ is a character and so $\fc(\eta\chi)  \leq  \max(\fc(\eta),\fc(\chi))$.
 Consequently, if $\pi$ is a principal series representation of $\GL_r(\Q_p)$ attached to the characters $\mu_1, ..., \mu_r$, we have
       $$ 
       \fc(\pi\otimes\chi)  =  \prod_j  \fc(\mu_j\chi)  \leq  \prod_j  \max(\fc(\mu_j),\fc(\chi)).
       $$
       
Now take $t=2$ or $t=3$.   As $ \eta$ is a discrete series it is either Steinberg $St$, or a twisted Steinberg $St(\nu)$ for a character $\nu$ of $\Q_p$, or a supercuspidal representation.
We then have the following situations:
\begin{itemize}
\item
If $\eta=St$  then $\fc(\eta)=p$ and $ \fc(\eta\otimes\chi)  =  \fc(\chi) = p^b$.
\item
If $\eta = St(\nu)$ then  $\fc(\eta)=\fc(\nu)$ and  $\fc(\eta\otimes\chi) \fc(St(\nu\chi)) = \fc(\nu\chi)$  in the case  $\nu\chi \ne 1$, while $\fc(\eta\otimes\chi) \fc(St(\nu\chi)) = p$ if $\nu\chi=1$.
\item
If $ \eta$ is supercuspidal (recall we assumed
$ t < p$), then %as by hypothesis, 
 $\eta$ is attached to a character $\lambda$ of a cyclic $t$-extension $K$ of $\Q_p$.  We have
                $\fc(\eta) = N(\fc(\lambda))d_K$, where $N$ denotes the norm from $K$ to $\Q_p$, and $d_K$ denotes the discriminant of $K/\Q_p$.
                Moreover,
                $\eta\otimes\chi$  is attached to the character $\lambda(\chi\circ N)$ of $K$, so that
               $$ \hskip.5in \fc(\eta\otimes\chi) = N(\fc(\lambda(\chi\circ N)))d_K  \leq  N\max(\fc(\lambda), \fc(\chi\circ N))d_K.$$
               
            There are really two types of supercuspidals $\eta$, depending on whether $K/\Q_p$ is unramified or ramified. In the former case, $d_K=1$ and 
            $p$ defines a uniformizer of $K$, so that (in this case) $\fc(\eta)= \fc(\lambda)$ and
                 $\fc(\eta\otimes\chi) =   \max (\fc(\lambda), p^{k/t})$ if  $\fc(\chi)=p^k$.
            Next suppose $K/\Q_p$ is ramified with $d_K=p^x$.  In this case, if  $\varpi$ is a uniformizer of $K$, we have $p\cO_K = \varpi^t \cO_K$ and $N(\varpi)=p$.
            Thus, if $\fc(\lambda)=\varpi^j$, then 
            $\fc(\eta)= p^{j+x}$,
            and
           $\fc(\eta\otimes\chi) =  p^{x+\max(j, b)}$.
\end{itemize} 
This completes our analysis.

%\newpage

\section{The Functional Equations and Their Properties}\label{sec:FEs}

\subsection{The Interchange Property %Actual Functional Equation
}\label{sec:interch}
\

The interchange property 
alluded to
in
the transfer from \eqref{step1} to \eqref{step2} is a very well-developed
component of the theory of Multiple Dirichlet Series.
The exact ``correction'' polynomials and their properties 
at unramified places
are detailed 
 in many places, including %e.g. 
 \cite{GoldfeldHoffstein1985,
BumpFriedbergHoffstein1996, DiaconuGoldfeldHoffstein2003, BumpFriedbergHoffstein2004,  ChintaFriedbergHoffstein2006, BrubakerBumpChintaFriedbergHoffstein2006}, to name a few. 
We assume some familiarity with these sources, while making the 
observation
 discussed below, that the correction polynomials can also be defined at ramified places.

\subsubsection{The Standard Approach}
\

First we recall the standard approach to the interchange property.
Let $S$ be a finite set of primes consisting of $2$ and the primes dividing the level $N$ of $\pi$. Let $M:=\prod_{p\in S}p$, and let 
$$
L^{S}(s,\pi\otimes\chi):=\prod_{p\notin S}L_{p}(s,\pi\otimes\chi)
$$
denote the twisted $L$-series with the places dividing $M$ removed.
Let $\ell_{1},\ell_{2}\mid M$ with $\ell_{1},\ell_{2}>0$ and $a_{1},a_{2}\in\{-1,1\}$.
One defines
\be\label{sfunct}
Z^{S}(s,w;\chi_{a_{2}\ell_{2}},\chi_{a_{1}\ell_{1}};\pi)
:=
\sum_{(d,M)=1}
{
L^{S}
(s,\pi\otimes\chi_{d_{0}}
\chi_{a_{1}\ell_{1}}
) 
\chi_{a_{2}\ell_{2}}(d_{0}) 
P^{(a_{1}\ell_{1})}_{d_{0},d_{1}}(s)
\over
%(d\ell_{1})
d^{w}
}
,
\ee
where the sum is over $d>0$ and we use the decomposition $d=d_{0}d_{1}^{2}$ with $d_{0}$ square-free.
This allows further character twists inside  the $L$-function by $\chi_{a_{1}\ell_{1}}$ and in the numerator by $\chi_{a_{2}\ell_{2}}$.    
One only considers such sums in which all the $L$-series of the numerator share a common gamma factor.  (A given sum can be, if necessary,  subdivided into several 
such sums.)

%\begin{remark}
%Notice that certain 
The ``correction'' polynomials $P(s)$
 %have been introduced.   These 
appearing in \eqref{sfunct} 
 are finite Dirichlet series that are uniquely determined by certain functional equations and limiting values.   They are introduced to, in essence, extend the level of a character from a square-free $d_0$ to $d_0d_1^2$ while retaining a similar functional equation in which the level $d_0$ is replaced by $d_0d_1^2$. Note this is \emph{not} the case if the primitive $\chi$ of conductor $d_0$ is simply replaced by the imprimitive character of conductor  $d_0d_1^2$.  

Take, for example, the simplest $\GL(1)$ case,  in which $L_{p}(s,\pi\otimes\chi)$ is the $p$-part of a quadratic Dirichlet $L$-series.  Here if $d_0 =q$ and $d_1 = p$ are distinct primes congruent to $1$ modulo $4$, and $d = qp^2$, then

$$P_{q,p}^{(1)}(s) = 1 - \chi_q(p)p^{-s} + p^{1-2s}.$$

Note that this has a functional equation $P_{q,p}^{(1)}(1-s) = p^{2s-1}P_{q,p}^{(1)}(s)$ and that by the original functional equation for $L(s,\chi_q)$, the product satisfies the particularly nice functional equation

$$
(qp^2/\pi)^{s/2}\Gamma(s/2)L(s,\chi_q)P_{q,p}^{(1)}(s) = L^*(s,\chi_q) =L^*(1-s,\chi_q)
.
$$
%\end{remark}

The 
now standard
 fact is that the polynomials $P^{(a_{1}\ell_{1})}_{d_{0},d_{1}}(s)$ can be defined in such a way that
the interchange 
\eqref{step2}
can be accomplished, with a new set of correction polynomials $Q(w)$ on the other side.

\begin{prop}\label{prop:InterStd}
There exists a choice of the  polynomials $P(s)$ and $Q(w)$ such that
%if in addition the following interchange of summation is valid for
%$s$ and $w$ having sufficiently large real part:
\be%a
%\nonumber
Z^{S}(s,w;\chi_{a_{2}\ell_{2}},\chi_{a_{1}\ell_{1}};\pi)
= \sum_{(n,M)=1} \frac{L^{S}(w,\tilde
\chi_{n_0}\chi_{a_2\ell_2})\chi_{a_1\ell_1}(n_0)
c(n_0n_1^2)Q_{n_0,n_1}^{(a_2\ell_2)}(w)
}{n^s}.
%\\
\label{Qswitch}
\ee%a
Here
the sum is over $n=n_{0}n_{1}^{2}$ where %with 
$n_{0}>0$ is
squarefree,
and
 $\tilde\chi_{n_0}$ denotes the quadratic character with
conductor $n_0$ defined by $\tilde
\chi_{n_0}(*) =
\left(\frac{*}{n_0}\right)$. 
(Recall $2|M$ so $(2,n_0)=1$.)
\end{prop}
The polynomials $Q_{n_0,n_1}^{(a_2\ell_2)}(w)$ have functional equation properties similar to those of $P^{(a_{1}\ell_{1})}_{d_{0},d_{1}}(s)$.

These correction polynomials have been explicitly worked out and written down  for every case 
considered 
here. The important point is that they exist and they are unique,
see e.g. (2.2) and (2.3) of \cite{BumpFriedbergHoffstein2004}.   
 Their exact form is cumbersome and not particularly illuminating for 
 our purposes,
 though their combinatorial properties have fascinating connections
 to statistical mechanics, crystal bases,  ice models, etc., cf. \cite{BrubakerBumpFriedberg2010}.
The main points to bear in mind are that 
\begin{enumerate}
\item For any $d_0,n_0$, $P^{(a_{1}\ell_{1})}_{d_{0},1}(s) = Q_{n_0,1}^{(a_2\ell_2)}(w)=1$.  That is, the correction polynomials are trivial when coefficients have square free indices.

\item The $P$ and $Q$ polynomials have simple functional equations that are compatible with the $L$-series by which %that 
they are multiplied. % by.
\item They permit interchanges such as that transforming \eqref{sfunct} into \eqref{Qswitch}.
% to take place.
\item The sizes of the $P$ and $Q$ polynomials are sufficiently small that for fixed $d_0, n_0$ the sums
$$
\sum_{d_1} \frac{P^{(a_{1}\ell_{1})}_{d_{0},d_{1}}(s)}{d_1^{2s}} \quad\text{and} \quad
\sum_{n_1}\frac{Q_{n_0,n_1}^{(a_2\ell_2)}(w)}{n_1^{2w}}
$$
converge absolutely for $\Re s >1/2$ and $\Re w >1/2$.
\end{enumerate}

\subsubsection{Ramified Correction Polynomials}\

We now make
explicit 
the aforementioned 
observation that 
it is not necessary to remove 
ramified $L$-parts
 to determine polynomials $P$ and $Q$ that satisfy the above four properties.  Take, for example the case of $\GL(1)$.  Here  $\pi = \chi$, a character of level $N$.   Taking, for simplicity, $a_1 = \ell_1 = a_2 = \ell_2 = 1$, 
 we replace \eqref{sfunct} 
 by
 \be\label{sfunctsimp}
Z(s,w;1,1;\pi)
:=
\sum_{d \ge1}
{
L(s,\chi \chi_{d_{0}}) 
P^{(1)}_{d_{0},d_{1}}(s)
\over
%(d\ell_{1})
d^{w}
}
.
\ee
Here  $L(s,\chi \chi_{d_{0}})$ denotes the primitive $L$-series attached to the character  
$\chi \chi_{d_{0}}$.  Thus if $(d_0,N)=1$ the conductor of the $L$-series is $Nd_0$.  In this case 
$P^{(1)}_{d_{0},d_{1}}(s)$ retains 
its standard
 definition, even when $d_1$ is divisible by primes dividing $N$.    
 If $(d_0,N) \ne1$,
  then one writes $Nd_0 = N_0d_3d_4^2$, where $N_0$ is the conductor of the product $\chi \chi_{d_{0}}$ and $d_3$ is square free.   In this case the numerator of the $d_0d_1^2$ term is  $L(s,\chi \chi_{d_{0}})P^{(1)}_{d_{3},d_{4}}(s)$.   For example, if $N=p$ is prime and $\chi = \chi_p$ is the quadratic character modulo $p$, then for any
$d_1$, and $d_0 = pd_0'$, the numerator becomes $L(s, \chi_{d_0'})P^{(1)}_{d_{0}',pd_{1}}(s)$.  The numerator has been constructed to have the correct functional equation, and all that needs to be verified is the interchange property, which is easily done.   In particular, in place of
\eqref{Qswitch}
we obtain
\be\label{Qswitchsimp}
Z(s,w;1,1;\pi)
=
\sum_{n \ge1}
{
L(w,\chi_{n_{0}})\chi(n)Q_{n_0,n_1}^{(1)}(w)
\over
%(d\ell_{1})
n^{s}
}
.
\ee

The reason why this interchange still holds, even when divisibility at ramified places is not restricted,  is that the $P$ and $Q$ polynomials exist because of a uniqueness principle that holds in $\GL(r)$ when $r = 1,2,3$.  See 
\cite[\S3.5]{ChintaFriedbergHoffstein2006} %Section 3.5 
for an exposition of this principle and a description of how the properties above are used to determine the  $P$ and $Q$ polynomials.  In $\GL(2)$ the $p$-part of an
$L$-series corresponding to a ramified prime can have one or zero Satake parameters.   If zero, the same discussion as in $\GL(1)$ above can be used to determine the $p$-parts of $P$ and 
$Q$ polynomials.   If one, then the the $p$-part is identical to the corresponding  parts of the polynomials in the $\GL(1)$ context, with the parameter $\alpha = \pm 1$ replacing $\chi(p)$.   Similarly, in $\GL(3)$, at ramified primes the $p$-parts of the $P,Q$ polynomials reduce to their $\GL(1)$ and $\GL(2)$ counterparts corresponding to $0,1$ or $2$ Euler factors.

In conclusion, we have the following (cf. Proposition \ref{prop:InterStd}).
\begin{prop}
For $\ell_{1},\ell_{2}\in\{1,2\}$ and $a_{1},a_{2}\in\{\pm1\}$, let
$$
Z(s,w;\chi_{a_{2}\ell_{2}}, \chi_{a_{1}\ell_{1}};\pi)
:=
\sum_{d\ge1}
{
L(s,\pi\otimes\chi_{d_{0}}\chi_{a_{1}\ell_{1}})
\chi_{a_{2}\ell_{2}}(d_{0})
P^{(a_{1}\ell_{1})}_{d_{0},d_{1}}(s)
\over
d^{w}
},
$$
where the $P$ polynomial is as described above. Then there is a $Q$ polynomial 
satisfying the above properties so that
$$
Z(s,w;\chi_{a_{2}\ell_{2}}, \chi_{a_{1}\ell_{1}};\pi)
=
\sum_{n\ge1}
{
L(w,\tilde\chi_{n_{0}}\chi_{a_{2}\ell_{2}})
\chi_{a_{1}\ell_{1}}(n_{0})
c(n_{0}n_{1}^{2})
Q^{(a_{2}\ell_{2})}_{n_{0},n_{1}}(w)
\over
n^{s}
}
.
$$
\end{prop}

%There is a rather large amount of notation involved in the interchange of summation formula above.   
As the presence 
of the correction polynomials is distracting,
we %define
fix the notation
\be\label{Pdef}
L(s,\pi\otimes\chi_d\chi_{a_1\ell_1}) :=    
L(s,\pi\otimes\chi_{d_0}\chi_{a_1\ell_1})
P_{d_0,d_1}^{(a_1\ell_1)}(s)
\ee
and
\be\label{Qdef}
 L(w,\tilde
 \chi_{n}\chi_{a_2\ell_2}) :=    
 L(w,\tilde
 \chi_{n_0}\chi_{a_2\ell_2})
Q_{n_0,n_1}^{(a_2\ell_2)}(w)
.
\ee
Thus the interchange in the order of summation above takes the slightly more reasonable form:
\bea
\label{newinterchange}
\hskip-.5in
Z(s,w;\chi_{a_{2}\ell_{2}},\chi_{a_{1}\ell_{1}};\pi)
&=&
\sum_{d\ge1} \frac{L(s,\pi\otimes\chi_{d}\chi_{a_1\ell_1})\chi_{a_2\ell_2}(d)}{d^w} 
\\
\nonumber
&=& \sum_{n\ge1} \frac{L(w,\tilde
\chi_{n}\chi_{a_2\ell_2})\chi_{a_1\ell_1}(n)
c(n)}{n^s}
.
\eea

\subsection {The meromorphic continuation}
\

Modulo the caveat in the previous subsection about 
ramified correction polynomials, it is now a standard matter
to meromorphically continue the series above. For the reader's convenience, 
we include a sketch of the argument, keeping careful track of the dependence on the level of each occurring transformation.

Proceeding as in \cite{DiaconuGoldfeldHoffstein2003}, collect the set of quadratic characters appearing above into
$$
\cM:=\{\chi_{a%_{1} 
\ell%_{1}
} : a%_{1}
=\pm 1 ,\ \ell=1,2\},
$$
and assemble %collect 
the corresponding double Dirichlet series into the $16$-dimensional vector %s
$$
\vbZ(s,w;
\pi)
=
\left\{
Z(s,w;\chi_{a_{2}\ell_{2}},\chi_{a_{1}\ell_{1}};\pi)
\right\}_{\chi_{a_{1}\ell_{1}}\in\cM\atop\chi_{a_{2}\ell_{2}}\in\cM}.
$$
Let $\ga$ and $\gb$ denote the involutions $\ga:(s,w)\mapsto(1-s,w+r(s-1/2))$ and $\gb:(s,w)\mapsto(s+w-1/2,1-w)$.
The double Dirichlet series $\vbZ(s,w;\pi)$ converges absolutely in the tube region where the real parts of $s$ and $w$ 
exceed
%are greater than 
$1$.   Also, by the analytic continuation and functional equation of $L(s,\pi\otimes\chi_{d}\chi_{a_1\ell_1})$, it converges for all $s$ as long as the real part of $w$ is sufficiently large.   Similarly, by the interchange in \eqref{newinterchange},  it converges for all $w$ as long as the real part of $s$ is sufficiently large.  Let $R_1$ denote the tube region which is the union of all such $s,w$.

  From \eqref{newinterchange} we see that there are potential polar lines at $s = 1$ and $w =1$.  If the orders of these poles are $p_1,p_2$ respectively then  
 $(1-s)^{p_1} (1-w)^{p_2}  \vbZ(s,w;\pi)$ has a holomorphic continuation to $R_1$.  
 
  We now write $Z(s,w;\chi_{a_{2}\ell_{2}},\chi_{a_{1}\ell_{1}};\pi)$ in its interchanged form
 \benn%a
Z(s,w;\chi_{a_{2}\ell_{2}}\psi,\chi_{a_{1}\ell_{1}};\pi)
= \sum_{n \ge 1} \frac{L(w,\tilde
\chi_{n}\chi_{a_2\ell_2})\chi_{a_1\ell_1}(n)
c(n)}{n^{s}}.
\eenn%a 

Let $G_{a_2\ell_2,n}(w)$ be the gamma function associated to the Dirichlet $L$-series in the numerator.    This will depend upon $a_2\ell_2$ and the congruence class of $n$ modulo 8.  
Multiplying and dividing each term by the appropriate 
%by 
$G_{a_2\ell_2 }(w)$ and applying the 
$\beta$ involution (effectively applying the $r=1$ instance of  \eqref{functeq}), we obtain
\bea
 \nonumber
 \hskip-.5in
       Z(s,w;\chi_{a_{2}\ell_{2}},\chi_{a_{1}\ell_{1}};\pi) %\nonumber\\
   & =&( \delta_1(a_2\ell_2))^{1/2-w} \sum_{n \ge 1}\frac{G_{a_2 \ell_2,n}(1-w)}{G_{a_2 \ell_2,n}(w)} 
 \\ 
\label{beta}
& &\times\ \frac{L(1-w,\tilde
\chi_{n}\chi_{a_2\ell_2})\chi_{a_1\ell_1}(n)
c(n)}{n^{s+w-1/2}}.
\eea 
Here $ \delta_1(a_2\ell_2)$ is the power of 2 associated to the discriminant of the corresponding quadratic field.
Applying $\gb$  to $R_1$ maps $R_1$ into a new tube region which intersects $R_1$.    Inspecting \eqref{beta} we see that two potential additional polar lines have been  added:
$w = 0$ and $s+w-1/2 =1$.  Canceling these lines we see that 
$$(1-s)^{p_1} (1-w)^{p_2} (w)^{p_2}(s+w-1/2)^{p_1} \vbZ(s,w;\pi)$$ 
has a holomorphic continuation to $R_1$.  

Thus the original function, after the polar lines are cancelled, is extended to a new function which is holomorphic in the region $R_2$, 
defined to be
%which is 
the union $R_1 \cup \gb (R_1)$.  
 Note that no new poles are introduced by the gamma functions in the numerator as their 
 products
 % of this 
 with the
corresponding  
$L$ series 
are
%in the numerator is 
analytic, except possibly at 0 and 1.

The catch is that the order of summation can not be changed in \eqref{beta} as the sum over $n$ has a varying gamma ratio coefficient, depending upon the congruence class of $n$ modulo 8.
This is easily remedied by replacing the congruence class modulo 8 conditions by linear combinations of sums over all $n$ twisted by characters modulo 8.  Thus for each $a_1\ell_1, a_2 \ell_2$, the right hand side of \eqref{beta} breaks up into a linear combination summed over $a$ modulo 8, with constant coefficients, of pieces of the form
\bea
\label{beta2}
&&
\hskip-.5in
 ( \delta_1(a_2\ell_2))^{1/2-w} \frac{G_{a_2 \ell_2,a}(1-w)}{G_{a_2 \ell_2,a}(w)} 
\sum_{n \ge 1}\frac{L(1-w,\tilde
\chi_{n}\chi_{a_2\ell_2})\chi_{a_1'\ell_1'}(n)
c(n)}{n^{s+w-1/2}}
\\
\nonumber 
 &=&( \delta_1(a_2\ell_2))^{1/2-w} \frac{G_{a_2 \ell_2,a}(1-w)}{G_{a_2 \ell_2,a}(w)} 
 Z(s+w-1/2,1-w;\chi_{a_{2}\ell_{2}},\chi_{a_{1}'\ell_{1}'};\pi). 
 \eea 
The interchange property can now be applied to each of the above pieces.
For a further discussion of this detail,
see \cite{DiaconuGoldfeldHoffstein2003},  line (4.23) in the arxiv version, and the discussion just preceding it. 

Now let
$G^{(a_1 \ell_1)}_\pi(s+w-1/2)$ be the gamma factor associated  to the $L$-series  
in the numerator of the 
reflected 
series in \eqref{beta2}.   This will be the common gamma factor for all $d \ge 1$.
 Multiplying and dividing by $G^{(a_1 \ell_1)}_\pi(s+w-1/2)$ and using the known functional equation of the $L$-series in the $s$ variable (see \eqref{functeq}), we effectively apply the involution $\ga$, obtaining for each $a_1\ell_1$,
\bea
\label{alpha}
&&
\hskip-.5in
       Z(s+w-1/2,1-w;\chi_{a_{2}\ell_{2}},\chi_{a_{1}\ell_{1}};\pi) 
       \\
 \nonumber
       & =&
   (N \delta_2(a_1\ell_1))^{1-s-w}\frac{G^{(a_1 \ell_1)}_\pi(3/2-s-w)}{G^{(a_1 \ell_1)}_\pi(s+w-1/2)} 
 \\ 
 \nonumber
 &&\times\ \ 
 Z(3/2-s-w,rs +(r-1)w+1-r;\chi_{a_{2}\ell_{2}}\psi,\chi_{a_{1}\ell_{1}};\pi).\nonumber
\eea 
Here $ \delta_2(a_1\ell_1)$ is again the power of 2 associated to the discriminant of the corresponding quadratic field.  Also we must be careful to choose regions for the parameters where the sum is absolutely convergent in both variables.  Notice that the character 
$\chi_{a_{2}\ell_{2}}$ has been multiplied by the character $\psi$ introduced by the application of the functional equation. 

We now have two other potential polar line at $s+w-1/2 =0$ and at $rs +(r-1)w+1-r = 1$.  Canceling these lines we see that 
\beann
\label{r=1}
&&(1-s)^{p_1} (1-w)^{p_2} (w)^{p_2}(s+w-1/2)^{p_1}(3/2 -s-w)^{p_1}\\
&&\times
\ \ 
(r-rs-(r-1)w)^{p_2} \vbZ(s,w;\pi)
\eeann
has a holomorphic continuation to $R_3$, the union $R_2 \cup \ga (R_2)$.  
Again no new poles are introduced by the gamma functions in the numerator 
since
 the product of this with the $L$-series in the numerator is analytic (except possibly at 0 and 1).

If $r=1$, the argument of $Z$ has been transformed from $(s,w)$ to $(3/2-s-w,s)$.  In this case the convex hull of $R_3$ is all of $\C^2$.  A  theorem of 
B\"ochner, see e.g.
%H\"ormander 
\cite[Thm 2.5.10]{Hormander1990}, %Theorem 2.5.10, 
then extends the domain of holomorphy of the function in \eqref{r=1} to this convex closure, namely $\C^2$.   
The original function 
$\vbZ(s,w;\pi)$ extends to a meromorphic function of $s,w$ in $\C^2$ with poles cancelled by 
$$
P_1(s,w) = s(1-s) w(1-w) (s+w-1/2)(3/2 -s-w)
,
$$
 and with polynomial growth in vertical strips determined, as usual, by the convexity principle.
A very detailed explanation of this process in the case $r=3$, along with illustrations,  is given in  \cite
[Proposition 4.11]
{DiaconuGoldfeldHoffstein2003}. 
%, Proposition 4.11. 

In the case $r=2$ an additional application of the $\beta$ involution is applied to 
$Z(3/2-s-w,w+2s-1;\chi_{a_{2}\ell_{2}}\psi,\chi_{a_{1}\ell_{1}};\pi)$ in \eqref{alpha}.  
Letting $N'$, with $N'|N$, denote the conductor of $\psi$, this leads to
\bea
\label{beta3}
       &&
 \hskip-.5in      Z(3/2-s-w,w+2s-1;\chi_{a_{2}\ell_{2}}\psi,\chi_{a_{1}\ell_{1}};\pi) 
 \\
 \nonumber
   & =&
   (N' \delta_3(a_2\ell_2\psi))^{3/2-2s-w} \sum_{n \ge 1}\frac{G_{a_2 \ell_2,n}(2-2s-w)}{G_{a_2 \ell_2,n}(w+2s-1)} 
 \\
 \nonumber
  &&\times\ \  \
  \frac{L(2-2s-w,
 \tilde \chi_{n}\chi_{a_2\ell_2}\overline{\psi})\chi_{a_1'\ell_1'}(n)\psi'(n)
c(n)\tau(\psi)}{n^{s}}.
\eea 
Here $\tau(\psi)$ is a (normalized to have absolute value 1) Gauss sum corresponding to $\psi$, $\psi'$ is a new character with conductor dividing $N$ arising from the factorization of the Gauss sum $\tau( \tilde \chi_{n}\chi_{a_2\ell_2}\psi)$  and $\chi_{a_1'\ell_1'}$ is a possibly new quadratic character modulo 8.  If $\psi^2 =1$, i.e if $\psi$ is trivial or quadratic, then $\psi'=1$.  %This is where 
It is here that
the case $\psi^2 =1$ begins to diverge from the case $\psi^2 \ne 1$.

Another sieving modulo 8 is now necessary to interchange the order of summation on the right hand side of \eqref{beta3}, transforming the right hand side into a linear combination of terms of the form
\bea\label{beta4}
      &&
      \hskip-.5in
      (N' \delta_3(a_2\ell_2\psi))^{3/2-2s-w} \sum_{n \ge 1}\frac{G_{a_2 \ell_2,n}(2-2s-w)}{G_{a_2 \ell_2,n}(w+2s-1)} 
 \\ &\times&Z(s,2-2s-w;\chi_{a_{2}\ell_{2}}\psi,\chi_{a_{1}\ell_{1}}\psi';\pi)\nonumber
\eea 

The convex hull of  $R_3 \cup \gb (R_3)$ is now, in the case $r=2$, all of $\C^2$.   Thus after canceling new potential polar lines we have the holomorphic continuation of 
$P_2(s,w) Z(s,2-2s-w;\chi_{a_{2}\ell_{2}},\chi_{a_{1}\ell_{1}};\pi)$ to all of $\C^2$, where
\bea
P_2(s,w)
& =& 
(s(1-s) (s+w-1/2)(3/2 -s-w))^{p_1} \nonumber 
\\ 
&&\times \ \
(w(1-w) 
(w+2s-1(2-2s-w))^{p_2}. \nonumber
\eea

In the case $r=3$, as when $r=2$, an additional application of the $\beta$ involution is applied to 
$Z(3/2-s-w,3s+2w-2;\chi_{a_{2}\ell_{2}}\psi,\chi_{a_{1}\ell_{1}};\pi)$ in \eqref{alpha}.  
This leads to
\bea
\nonumber
       &&
       \hskip-.5in
       Z(3/2-s-w,3s+2w-2;\chi_{a_{2}\ell_{2}}\psi,\chi_{a_{1}\ell_{1}};\pi)
        \\
\label{beta5}   & =&(N' \delta_3(a_2\ell_2\psi))^{5/2-3s-2w} \sum_{n \ge 1}\frac{G_{a_2 \ell_2,n}(3-3s-2w)}{G_{a_2 \ell_2,n}(3s+2w-2)} 
 \\
  &&\times \ \ \ \frac{L(3-3s-2w,
 \tilde \chi_{n}\chi_{a_2\ell_2}\overline{\psi})\chi_{a_1'\ell_1'}(n)\psi'(n)
c(n)\tau(\psi)}{n^{w+2s-1}}.\nonumber
\eea 
Here $\tau(\psi)$,$\psi' $, $N'$   and $\chi_{a_1'\ell_1'}$ are as in the case $r=2$.   As in that case, $\psi' =1$ whenever $\psi$ is trivial or quadratic.

Another sieving modulo 8 is now necessary, as when $r=2$,  to interchange the order of summation on the right hand side of \eqref{beta3}, transforming the right hand side into a linear combination of terms of the form
\bea\label{beta6}
      &&
      \hskip-1in
      (N' \delta_3(a_2\ell_2\psi))^{5/2-3s-2w} \sum_{n \ge 1}\frac{G_{a_2 \ell_2,n}(3-3s-2w)}{G_{a_2 \ell_2,n}(3s+2w-2)} 
 \\ &\times&Z(w+2s-1,3-3s-2w;\chi_{a_{2}\ell_{2}}\psi,\chi_{a_{1}\ell_{1}}\psi';\pi)
.
 \nonumber
\eea 

We now apply the $\ga$ involution to $Z(w+2s-1,3-3s-2w;\chi_{a_{2}\ell_{2}}\psi,\chi_{a_{1}\ell_{1}}\psi';\pi)$, obtaining
\bea\label{alpha2}
       &&
       \hskip-.5in
       Z(w+2s-1,3-3s-2w;\chi_{a_{2}\ell_{2}}\psi,\chi_{a_{1}\ell_{1}}\psi';\pi) \nonumber \\
   & =&(N'' \delta_4(a_1\ell_1\psi'))^{3/2-2s-w}\frac{G^{(a_1 \ell_1)}_\pi(2-2s-w)}{G^{(a_1 \ell_1)}_\pi(w+2s-1)} 
 \\ 
 &&\times \ \ \ 
 Z(2-2s-w,3s +w-3/2;\chi_{a_{2}\ell_{2}}\psi'',\chi_{a_{1}\ell_{1}}\overline{\psi'};\pi).\nonumber
\eea 
Here $N''$ is the conductor of $\pi \otimes \psi'$, $\psi''$ is the new functional equation character twist of $\pi \otimes \psi'$  multiplied by $\psi$, and $\delta_4(a_1\ell_1\psi')$ keeps track of the extra powers of $2$ introduced by the twisting.   We refer to the analysis of 
\S\ref{sec:Din}
for the computation of $N''$ in the general case.  However, if $\psi^2 = 1$ then $\psi' =1$ and $\psi'' = \psi^2 =1$. In this case  $\psi''$ is an imprimitive identity character modulo  $N'$ and $N'' = N$.

We now need to apply the involution $\gb$ one more time, to 
$Z(2-2s-w,3s +w-3/2;\chi_{a_{2}\ell_{2}}\psi'',\chi_{a_{1}\ell_{1}}\overline{\psi'};\pi)$.   
   As a consequence the $\GL(1)$ $L$-series in the $3s +w-3/2$ variable will be missing Euler factors at primes dividing $N'$, the conductor of $\psi$.  For this reason we write 
 \bea
&&
\hskip -.5in
Z(2-2s-w,3s +w-3/2;\chi_{a_{2}\ell_{2}}\psi'',\chi_{a_{1}\ell_{1}}\overline{\psi'};\pi)\nonumber \\
&=& \sum_{n \ge 1} \frac{L(3s +w-3/2,
\tilde\chi_{n}\chi_{a_2\ell_2}\psi_0'')
\chi_{a_1\ell_1}(n)
c(n)}{n^{2-2s-w}} \nonumber\\
&&\times \ \ \ \prod_{p|N'}\big(1-\tilde\chi_{n}\chi_{a_2\ell_2}\psi_0''(p)p^{3/2-3s-w}\big),
\eea 
where $\psi_0''$ is the corresponding primitive character.

Applying $\gb$, and then sieving modulo 8, we end up with a linear combination of terms of the form

\bea\label{beta6}
      &&
      \hskip-1in(N''' \delta_5(a_2\ell_2\psi_0''))^{2-3s-w}(N'''')^{3/2-3s-w} \sum_{n \ge 1}\frac{G_{a_2 \ell_2,n}(5/2-3s-w)}{G_{a_2 \ell_2,n}(w+3s-3/2)} \nonumber
 \\ &\times&Z(s,5/2-3s-w;\chi_{a_{2}\ell_{2}}\overline{\psi_0''},\chi_{a_{1}\ell_{1}}\psi';\pi)
 .
 \eea 
Here $N'''$ is the conductor of the $\GL(1)$ $L$-series of the character $\psi'$ (so $N''' =1$ if 
$\psi'=1$) and $N''''$ is the conductor of the imprimitive part of $\psi''$.

We  have now continued to a region whose convex closure is $\C^2$.   Canceling the potential polar lines, we finally have, in the case $r=3$, a holomorphic continuation of 
$P_3(s,w) Z(s,w;\chi_{a_{2}\ell_{2}},\chi_{a_{1}\ell_{1}};\pi)$ to all of $\C^2$, where
\bea
P_3(s,w) &=& (s(1-s) (s+w-1/2)(3/2 -s-w)(2s+w -1))^{p_1} \nonumber 
\\ 
&&\times \ \ \ 
(2-2s-w)^{p_1}(w(1-w)
(3s+2w-2)(3-3s-2w))^{p_2}. \nonumber
\eea

The above calculations are all in the literature, but included here for the reader's convenience, and for the exact dependence on the various occurrences of $N$, $N'$, \dots, $N''''$.

\subsection{Collecting the analytic information}
\

Setting $s=1/2$ in the information gathered above enables us to give a precise description of the analytic behavior of the functions $Z(1/2,w;\chi_{a_{2}\ell_{2}},\chi_{a_{1}\ell_{1}};\pi)$ as follows
\begin{proposition}
\label{prop:convexbounds}
The function  
$$
Z(1/2,w;\chi_{a_{2}\ell_{2}},\chi_{a_{1}\ell_{1}};\pi)
=
\sum_{d \ge1} \frac{L(1/2,\pi\otimes\chi_{d}\chi_{a_1\ell_1})\chi_{a_2\ell_2}(d)}{d^w} 
$$
converges absolutely when $\Re w >1$ and has a meromorphic continuation to $\C$.  
It has poles at
\begin{enumerate}
\item
$w=1$,  when $r=1$,
\item
$w=1$,  when $r=2$,
\item
$w=1$ and $w=3/4$ when $r=3$. 
\end{enumerate}
Note that the function ``completed'' with Gamma factors (see \eqref{funct1new} -- \eqref{funct3new}) also has poles at $w=0$, and on $\GL(3)$ at $w=1/4$.
Away from these poles it has polynomial growth in vertical strips.
   
Finally, for any $a_1\ell_1,a_2\ell_2\in\{\pm1,\pm2\}$, it satisfies the following functional equations.   
The functions $C_{1}$, $C_{2}$, and $C_{3}$ occurring below are absolute
constants depending on values in $\{\pm1,\pm2\}$.
When  $r=1$,
\bea
\label{funct1}
&&
\hskip-.5in Z(1/2,w;\chi_{a_{2}\ell_{2}},\chi_{a_{1}\ell_{1}};\pi)\\
&=& \sum_{a_{1}'\ell_{1}',a_{2}'\ell_{2}'}C_1(a_1\ell_1, a_2\ell_2,a_1'\ell_1', a_2'\ell_2')( \delta_1(a_2'\ell_2'))^{1/2-w} \frac{G_{a_2' \ell_2'}(1-w)}{G_{a_2' \ell_2'}(w)} \nonumber
\\
&&\times \ \ \ (N \delta_2(a_1'\ell_1'))^{1/2-w}\frac{G^{(a_1' \ell_1')}_\pi(1-w)}{G^{(a_1' \ell_1')}_\pi
(w)} 
%\nonumber\\ &\times&
Z(1-w,s ;\chi_{a_{2}'\ell_{2}'}\psi,\chi_{a_{1}'\ell_{1}'};\pi)  
\nonumber
.
\eea
When  $r=2$,
\bea\label{funct2}
&&
\hskip-.5in
Z(1/2,w;\chi_{a_{2}\ell_{2}},\chi_{a_{1}\ell_{1}};\pi)\\
&=& \sum_{a_{1}'\ell_{1}',a_{2}'\ell_{2}'}C_2(a_1\ell_1, a_2\ell_2,a_1'\ell_1', a_2'\ell_2')( \delta_1(a_2'\ell_2'))^{1/2-w} \frac{G_{a_2' \ell_2'}(1-w)}{G_{a_2' \ell_2'}(w)} \nonumber \\
&&\times\ \ (N \delta_2(a_1'\ell_1'))^{1/2-w}\frac{G^{(a_1' \ell_1')}_\pi(1-w)}{G^{(a_1' \ell_1')}_\pi
(w) 
}
%\nonumber\\ 
% &\times&
 (N' \delta_3(a_2'\ell_2'\psi))^{1/2-w} 
 \frac{G_{a_2' \ell_2'}(1-w)}{G_{a_2' \ell_2'}(w)} 
   \nonumber \\ 
&&\times \ \ 
Z(1/2,1-w;\chi_{a_{2}'\ell_{2}'}\psi,\chi_{a_{1}'\ell_{1}'}\psi';\pi) \nonumber
.
\eea
When  $r=3$, 
\bea\label{funct3}
&&
\hskip-.5in
Z(1/2,w;\chi_{a_{2}\ell_{2}},\chi_{a_{1}\ell_{1}};\pi) \\
&=& \sum_{a_{1}'\ell_{1}',a_{2}'\ell_{2}'}C_3(a_1\ell_1, a_2\ell_2,a_1'\ell_1', a_2'\ell_2')( \delta_1(a_2'\ell_2'))^{1/2-w} \frac{G_{a_2' \ell_2'}(1-w)}{G_{a_2' \ell_2'}(w)}  \nonumber\\
&&\times  \ \ (N \delta_2(a_1'\ell_1'))^{1/2-w}\frac{G^{(a_1' \ell_1')}_\pi(1-w)}{G^{(a_1' \ell_1')}_\pi
(w)} 
%\nonumber\\ 
 %&\times&
 (N' \delta_3(a_2'\ell_2'\psi))^{1-2w} %\sum_{n \ge 1}
 \frac{G_{a_2' \ell_2'}(3/2-2w)}{G_{a_2' \ell_2',}(2w-1/2)} 
 \nonumber 
 \\ 
 &&\times \ \ 
 (N'' \delta_4(a_1'\ell_1'\psi'))^{1/2-w}\frac{G^{(a_1' \ell_1')}_\pi(1-w)}{G^{(a_1' \ell_1')}_\pi(w)}  
 (N''' \delta_5(a_2'\ell_2'\psi_0''))^{1/2-w}(N'''')^{-w} 
 \nonumber\\ 
 &&\times \ \ \
\frac{G_{a_2' \ell_2'}(1-w)}{G_{a_2' \ell_2'}(w)} 
%  \nonumber
%\\ &\times&
Z(1/2,1-w;\chi_{a_{2}'\ell_{2}'}\overline{\psi_0''},\chi_{a_{1}'\ell_{1}'}\psi';\pi) 
.
\nonumber
\eea

For $w = -\epsilon + it$, with $\epsilon >0$, lines \eqref{funct1} -- \eqref{funct3} relate 
$Z(1/2,-\epsilon + it;\chi_{a_{2}\ell_{2}},\chi_{a_{1}\ell_{1}};\pi)$ to
$Z(1/2,1+\epsilon - it;\chi_{a_{2}\ell_{2}},\chi_{a_{1}\ell_{1}};\pi)$.  Suppose that
$$
Z(1/2,1+\epsilon - it;\chi_{a_{2}\ell_{2}},\chi_{a_{1}\ell_{1}};\pi) \ll_{t,\fq} N^{\alpha_r},
$$
with
some $\ga\ge0$, and
 the implied constant 
 depending at most
  polynomially
  on  $t$.   
  %Here $\alpha_r \ge 0$ can be set to $0$ if the Lindel\"of Hypothesis is true in the $N$ aspect, or to the best known upper bound.
Then by convexity it follows that for $r=1$,
\be\label{eq:convR1}
Z(1/2,1/2+ it;\chi_{a_{2}\ell_{2}},\chi_{a_{1}\ell_{1}};\pi) \ll N^{1/4 + \alpha_1+\epsilon},
\ee
for $r=2$,
\be\label{eq:convR2}
Z(1/2,1/2+ it;\chi_{a_{2}\ell_{2}},\chi_{a_{1}\ell_{1}};\pi) \ll (NN')^{1/4+ \alpha_2 + \epsilon},
\ee
and for $r=3$,
\be\label{eq:convR3}
Z(1/2,1/2+it;\chi_{a_{2}\ell_{2}},\chi_{a_{1}\ell_{1}};\pi) \ll (N(N')^2N'' N''')^{1/4 + \alpha_3+ \epsilon},
\ee
with the implied constant including a polynomial power of $t$.

Here $N$ is the conductor of $\pi$, $N'$ is the conductor of $\psi$,
$N''$ is the conductor of $\pi\otimes \psi'$, and $N'''$ is discussed just after \eqref{beta6}.
  If $\psi^2 = 1$ then $N'' = N$ and $N''' =1$.

\end{proposition}
The previous proposition is in a useful form for measuring the precise growth in vertical strips via an application of Stirling's formula to the ratios of gamma factors.     After multiplying by the gamma factors in the denominator, another more symmetrical formulation of the functional equation can be found.  This is the form that we will use in the proof of Theorem \ref{thm0}.

The gamma factors $G_{a_2' \ell_2'}(w),G_{a_2' \ell_2',n}(w)$ are all of the form
 $$
G_{+}(w) = \pi^{-w/2}\Gamma(w/2)\,\, \text{or} \,\, G_{-}(w) =(2\pi)^{-w/2}\Gamma((w+1)/2).
 $$
  Similarly, the gamma factors $G^{(a_1' \ell_1')}_\pi(w)$ are all of the form 
 $$
G_{+,\pi}(w) =\prod_{i=1}^r \pi^{-w/2}\Gamma((w+\kappa_i)/2)\,\, \text{or} \,\, G_{-,\pi}(w) =\prod_{i=1}^r (2\pi)^{-w/2}\Gamma((w+\kappa_i +1)/2).
 $$
Consider first the case $r=1$.   After multiplying by  $G_{+}(w)G_{+,\pi}(w)$, each piece in the summation on the right hand side of \eqref{funct1} takes the form a constant times
 $$
(N 2^a \pi^b )^{1/2-w}\rho (w)\rho_{\pi}(w)G_{+}(1-w)G_{+,\pi}(1-w)
Z(1-w,s ;\chi_{a_{2}'\ell_{2}'}\psi,\chi_{a_{1}'\ell_{1}'};\pi),
$$
where 
$$
\rho(w) = 1\,\, \text{or} \,\,G_{+}(w)/G_{-}(w)
$$ 
and
$$
\rho_{\pi} (w) =1 \,\, \text{or} \,\, G_{+,\pi}(w)/G_{-,\pi}(w).
$$
 Consequently \eqref{funct1} can be rewritten as
 \bea\label{funct1new}
&&
\hskip-.5in
G_{+}(w)G_{+,\pi}(w)Z(1/2,w;\chi_{a_{2}\ell_{2}},\chi_{a_{1}\ell_{1}};\pi) \\
&=& N^{1/2-w}\sum_{a_{1}'\ell_{1}',a_{2}'\ell_{2}'}( 2^a \pi^b )^{1/2-w}C_1\
\rho (w)\rho_{\pi}(w) \nonumber
\\ 
&&\times \ \ 
G_{+}(1-w)G_{+,\pi}(1-w)Z(1-w,s ;\chi_{a_{2}'\ell_{2}'}\psi,\chi_{a_{1}'\ell_{1}'};\pi).\nonumber
\eea
Here each $a,b,C_1,\rho, \rho_\pi$ is a function of $a_1\ell_1, a_2\ell_2$ and $a_1'\ell_1', a_2'\ell_2'$.

Similarly for $r=2$,  \eqref{funct2} can be rewritten as
 \bea\label{funct2new}
&&
\hskip-.5in
G_{+}(w)^2G_{+,\pi}(w)Z(1/2,w;\chi_{a_{2}\ell_{2}},\chi_{a_{1}\ell_{1}};\pi) \\
&=& (NN')^{1/2-w}\sum_{a_{1}'\ell_{1}',a_{2}'\ell_{2}'}( 2^a \pi^b )^{1/2-w}C_2\
\rho (w)^2\rho_{\pi}(w)\nonumber
\\ 
&&
\times \ \ \ 
G_{+}(1-w)^2G_{+,\pi}(1-w)Z(s,1-w ;\chi_{a_{2}'\ell_{2}'}\psi,\chi_{a_{1}'\ell_{1}'};\pi)
\nonumber,\eea
and for $r=3$,
\eqref{funct3} can be rewritten as
 \bea\label{funct3new}
&&
\hskip-.5in
G_{+}(w)^2G_{+}(2w-1/2)G_{+,\pi}(w)^2Z(1/2,w;\chi_{a_{2}\ell_{2}},\chi_{a_{1}\ell_{1}};\pi) \\
&=&  (NN'N''N''')^{1/2-w}(N'''')^{-w}\sum_{a_{1}'\ell_{1}',a_{2}'\ell_{2}'}( 2^a \pi^b )^{1/2-w}C_3\
\rho (w)^2 \rho(2w-1/2)\rho_{\pi}(w)^2 \nonumber
\\ 
&&\times\ \ \
G_{+}(1-w)^2G_{+}(3/2-2w)G_{+,\pi}(1-w)^2Z(s,1-w ;\chi_{a_{2}'\ell_{2}'}\psi,\chi_{a_{1}'\ell_{1}'};\pi).\nonumber
\eea

%\begin {remark}
\subsection{Residual Behavior}\label{sec:Siegel}\

The nature of the residues at the simple poles and leading coefficients of the Laurent expansion at higher order poles is rather complicated and varied,
and clearly necessary for our application of Tauberian arguments. Fortunately, the literature already contains sufficient information in each case occurring here.
Below, we collect the relevant facts.

For $s$ in a neighborhood of $1/2$ but $s\neq1/2$, the double Dirichlet series 
$Z(s,w;\chi_{a_{2}\ell_{2}},\chi_{a_{1}\ell_{1}};\pi)$ is holomorphic at $w=1$ as long as $a_{2}\ell_{2}\neq1$.  There are several polar lines intersecting the point $(1/2,1)$.  
The combination of these polar lines can create multiple poles or eliminate the pole at $(1/2,1)$, depending upon the number theoretic nature of $\pi$.

For example  in the case $r=2$ and $\pi$ a cuspidal newform these lines are $w=1$ and $w+s-1/2 =1$.  Each pole is simple, with residues a non-zero multiple of $L(2s,\pi,\sym^{2})$ and $L(1-2s,\pi,\sym^{2})$ respectively.  The residue of 
$Z(1/2,w;1,\chi_{a_{1}\ell_{1}};\pi)$ at $w=1$  is the sum of the limits of these two residues.   This sum is easily computed via 
%an easy calculation involving 
the global root number, and all cases it is zero if and only if the root number is $-1$.   This is the basis of non-vanishing theorems for families of quadratic $L$-series: in any case where the root number is not  $-1$ for all twists, the two residue terms do not cancel, forcing a pole which in turn forces 
the nonvanishing of
infinitely many quadratic twists of the relevant $L$-series at $1/2$. % to not vanish.
For a full discussion of this in the case $r=2$ see line (1.4) of \cite{FriedbergHoffstein1995} and the analysis following it.

In the case $r=3$ and $\pi$ cuspidal, if $a_2\ell_2=1$  and $s$ is in a neighborhood of $1/2$ but $s \ne 1/2$, 
then
the double Dirichlet series $Z(s,w;1,\chi_{a_{1}\ell_{1}};\pi)$ has a simple pole at $w=1$ of residue a constant multiple of 
$$
L(2s,\pi,\sym^2)\zeta(6s-1).
$$
%, where the $L$ series is the symmetric square $L$-series of $\pi$.   
If $\pi$ is itself the adjoint square lift of a cuspidal newform on $\GL(2)$, then 
%this $L$-series 
$L(2s,\pi,\sym^{2})$
has a pole at $s=1/2$, and hence $Z(1/2,w;1,\chi_{a_{1}\ell_{1}};\pi)$ has a double pole at $w=1$.  See the discussion after 
 Proposition 3.7 in \cite{BumpFriedbergHoffstein2004} for a full analysis of this case.
%\end{remark}

In all cases of interest the existence or non-existence of a pole of 
$Z(1/2,w;1,\chi_{a_{1}\ell_{1}};\pi)$ at $w=1$ can be checked.   In $\GL(1)$, $\GL(2)$ it has been verified that there is no pole at $w=1$ if and only if the root number of every twist of $\pi$ by a quadratic character is $-1$.   On $\GL(3)$, if $\pi$ is the adjoint square lift of a cuspidal newform on $\GL(2)$ there is always a pole of order 2.   For  generic $\pi$ on $\GL(3)$ the corresponding property can be easily checked.   Consequently we state the following %Proposition.
\begin{prop}
\label{prop:residue}
Suppose that $Z(1/2,w;1,\chi_{a_{1}\ell_{1}};\pi)$ has a pole of order at least 1 at $w=1$.   Denoting the coefficient of the leading term in the Laurant expression as $\kappa$,
there exists  some $A>0$ such that
\be\label{eq:ResBnd}
\kappa
\gg
(\fq N)^{-A}.
\ee
If the $L$-series $L(s, \pi, \sym^2)$ does not have a Siegel zero (see below for the definition) then the stronger result
\be\label{eq:ResBndSiegel}
\kappa
\gg
(\fq N)^{-\epsilon}
\ee
is true.
\end{prop}
\pf
If the pole is simple, then as discussed above, and also in the case $r=1$, $\kappa$ is equal to a non-zero multiple of $L(1, \pi, \sym^2)$.   In the case of a multiple pole, the relevant value is a lower rank $L$-series evaluated at 1.  It is well known that 
$$
\kappa \gg N^{-\epsilon} (1-\beta)
$$
for any $\beta$ satisfying
$$
1-(1- \log (\fq N))^{-1} <\beta <1
$$ 
and 
$L(s, \pi, \sym^2) \ne 0$ for  $s$ real with $\beta < s < 1$.
Such a $\beta$ is, if it exists, known as a Siegel zero.
Ordinarily, 
finding a lower bound for $1-\gb$ 
would be a subtle matter, and to obtain more refined results, it must be addressed.
But in this case,
 all we need is a polynomial lower bound,
 and hence even the weakest results in the cases $r=1, 2$ suffice.
Furthermore, the uniform bound in \cite{Brumley2006} gives the $\GL(3)$ claim, as the symmetric square $L$-function is a factor of the Rankin-Selberg convolution $L(2s,\pi\times\pi)$ and the exterior $L$-series factor is easily estimated from above.
\epf

\section{Proofs of 
Theorems \ref{thm0} and \ref{thm1}}
\label{sec:proof}

Let $\pi_r$ be an automorphic representation on $\GL(r)$, $r=1,2,3$ of level $N$  and let $\tilde \pi_r$ be its contragredient.   %Let us s
Suppose that there exists at least one quadratic character 
$\chi$ such that the root number of $\pi_r \otimes \chi$ is not equal to $-1$.

Let $G_r(w)$ denote the gamma factors %associated with 
which appear
in lines \eqref{funct1new} --  \eqref{funct3new}.   Thus 
$$G_1(w) = G_+(w)G_{+,\pi}(w),$$ 
$$G_2(w) = G_+(w)^2G_{+,\pi}(w)$$
 and 
$$G_3(w) = G_+(w)^2G_+(2w-1/2)G_{+,\pi}(w)^2.$$

For $X\ge1$, apply an inverse Mellin transform, obtaining
\be\label{eq:Iis}
\cI={1\over 2\pi i} \int_{(2)}G_r(w)  Z(1/2,w;1,1;\pi_r) X^{w} dw = \sum_{d}L(1/2,\pi\otimes \chi_d )  
V\left({d \over X}\right),
\ee
where 
$$
V(y):=
{1\over 2\pi i} \int_{(2)}G_r(w) y^{-w}dw
$$ 
has arbitrary polynomial decay 
for $y\gg1$.

Move 
the line of integration 
in \eqref{eq:Iis}
 to $\Re w = -1$,
 passing through poles  at $w=1$, $w=0$, and on $\GL(3)$, potential poles at  $w =3/4$ and $w = 1/4$.  Depending on $r$, apply the functional equations \eqref{funct1new} --  \eqref{funct3new}, and make the change of variables $w\mapsto 1-w$. 
 The moved integral reflects into the region of absolute convergence, and breaks into a linear combination of $L$-series $L(1/2,\pi\otimes \chi_d )$ times new damping functions
$\tilde V(dX/N^{2\theta_r})$.
The power $\gt_{r}$ of $N$ in the above is determined by the factors $N$, $N'$, \dots, $N''''$ in equations \eqref{funct1new} -- \eqref{funct3new}. For simplicity,
assume that  $\psi$  has conductor $N$, and
 moreover that if $r=3$, then $\psi$
 is  quadratic.
Then 
the power of $N$ is
$\gt_1 = 1/2,\ \gt_2 = 1$, and $\gt_3 = 2$.    
 
 Let $\gk_{1}$, $\gk_{3/4}$, $\gk_{1/4}$, and $\gk_{0}$ denote the residual contributions at $w=1$, $3/4$, $1/4$, and $0$, respectively.
By Proposition~\ref{prop:residue}, the leading term is of the order $\gg  XN^{-A}$  with a possible lower order residual terms $X^{3/4} N^{B}$, $X^{1/4}N^{C}$ and $N^{D}$ for some $A, B, C,$ and $D$. Let $\cR$ denote the residual contribution. By perturbing $X$ near 
$N^{\gt_{r}}$, we can ensure that  
$$
\cR\gg N^{-A}
$$
for some fixed $A>0$. Fix this value of $X$.   If there is no Siegel zero, then in fact 
$$
\cR\gg N^{-\epsilon}X.
$$

Assume now that for all $|d|\ll N^{\gt_{r}+\vep}$,  the twisted $L$-series $L(1/2,\pi\otimes\chi_{d})$  all vanish. 
 As $V$ and $\tilde V$  have arbitrary polynomial decay, 
 one obtains the desired
  contradiction, since the residual terms are $\gg N^{-A}$ for some $A$.
This contradiction
 completes the proof of Theorem \ref{thm0}.

For
Theorem \ref{thm1} in the eigenvalue aspect,  one can easily apply Stirling's formula to the gamma factors corresponding to the level and 
complete the proof 
using the same argument.

\begin{rmk}\label{rmk:SubConv}

Were one to attempt an improvement along the lines of \S\ref{sec:conv},
one could proceed as follows. Instead of \eqref{eq:Iis}, let $h(y)$ denote a smooth function of compact support in $[1,2]$, say, and let $H(w)$ denote its Mellin transform, $H(w)=\int_{0}^{\infty}h(y) y^{w} dy/y.$ The latter has arbitrary polynomial decay in vertical strips if $\Re(w)>0$.
Take an integral like \eqref{eq:Iis}, except with $H$ replacing the Gamma factors:
$$
\cI:=
\frac1{2\pi i} \int _{(2)}H(w) Z(1/2,w;1,1;\pi)X^{w}dw
=
\sum_{d\ge1}
L(1/2,\pi\otimes\chi_{d})h(d/X).
$$
Now
move the contour back to $\Re(w)=1/2$, passing through the poles  at $w=1$ (and possibly at $w=3/4$ on $\GL(3)$). 
Let $\kappa$ denote the residual contribution from the pole at $w=1$, $\kappa'$ the potential contribution from $w=3/4$.

Assuming there is no Siegel zero (which has been established for our purposes in all but the classical case), then again 
$$
\gk\gg N^{-\epsilon}.
$$

Then inputting the supposed bound \eqref{eq:conj} and estimating away the residual contribution gives
$$
 \sum_{d}L(1/2,\pi\otimes \chi_d ) \ h\left({d \over X}\right)\gg N^{-\vep} X 
 %+ \kappa' X^{3/4}
 + \cO \left(X^{1/2} N^{\epsilon}\right).
$$
Take $X=N^{4\vep}$.
Assuming that $L(1/2,\pi\otimes \chi_d ) =0$ for all $d \le 2 X$, 
one obtains the desired contradiction.
This completes our analysis.
\end{rmk}

\bibliographystyle{alpha}

\bibliography{AKbibliog}

\end{document}